\documentclass[a4paper, 12pt]{article}

\usepackage{graphicx}
\usepackage{url}
\usepackage{amsbsy,amssymb} 
\usepackage{amsmath}
\usepackage{abstract}
\usepackage{hyperref}
\usepackage{color}

\usepackage{amsmath,amsthm}

\newtheorem{Proposition}{Proposition}

\newtheorem{Corollary}{Corollary}

\newtheorem*{Proof}{Proof}

\begin{document}

{\LARGE\centering{\bf{Perturbed Lane-Emden equations as a boundary value problem with singular endpoints}}}

\begin{center}
\sf{Rados\l aw Antoni Kycia$^{1,2,a}$}
\end{center}

\medskip
\small{
\centerline{$^{1}$Masaryk Univeristy}
\centerline{Department of Mathematics and Statistics}
\centerline{Kotl\'{a}\v{r}sk\'{a} 267/2, 611 37 Brno, The Czech Republic}
\centerline{\\}
\centerline{$^{2}$Cracow University of Technology}
\centerline{Faculty of Physics, Mathematics and Computer Science}
\centerline{Warszawska 24, Krak\'ow, 31-155, Poland}
\centerline{\\}

\centerline{$^{a}${\tt
kycia.radoslaw@gmail.com}}
}

\begin{abstract}
\noindent
The paper presents the solution for the existence of analytic solutions for some generalized Lane-Emden (LE) equation. Such solutions exists on the unit interval, which endpoints are singularities of the proposed perturbed LE equation. The solution has many possible applications and one of the examples was provided.
\end{abstract}
Keywords: Generalized Lane-Emden equations, analytic solutions, singularities;  \\
Mathematical Subject Classification: 34A34, 34M35;\\

\section{Introduction}

The Lane-Emden equation (LE) \cite{Lane, EmdenGaskugeln, ChandrasekarBook} is one of the simplest nonlinear equations appearing in Astrophysics. It describes balance in gaseous medium described by polytropic equation, between qravity that tries to squeeze it and internal pressure that try to prevent such collapse.

It is used in Astrophysics as a first approximation to the more complicated star structure models \cite{ChandrasekarBook, Lane, EmdenGaskugeln, Hunter_LaneEmden} and recently for describing molecular cloud cores \cite{MolecularCloudCores}. Some exact solutions \cite{Mach_n5} are known and are useful in solving (matching) problems in which interior is described by LE equation and exterior area by some other equation. Such exact solutions appears also for higher space dimensions \cite{KyciaFilipukCriticalHigherDimensions} and can be used to representation of symmetric Riemannian manifold with constant curvature as conformal flat model \cite{MachCriticalLEHigherDimension}.

However, apart of applications, the LE equation is the simplest nonlinear equation that arises from radial part of the Laplace operator in flat space with nonlinear term. Therefore it serves also as a building block for more complicated nonlinear differential equations, e.g., nonlinear wave equation with power-type nonlinearity \cite{metoda_dim=3, wyzsze_wymiary}. There are various kinds of generalization for this equations described, e.g., in \cite{GeneralizedLE1, GeneralizedLE2, ISAAC, KyciaFilipuk_LandEmden}. These solutions poses interesting structure of movable singularities (singularities of analytic solutions) \cite{KyciaFilipuk_LandEmden, KyciaFilipuk_LandEmden, wyzsze_wymiary_movable_singularities} and some flavour of the Painlev\'{e} analysis for them was provided in \cite{GeneralizedLE2}.

The main aim of this paper is to perturb generalized the Lane-Emden equation \cite{Davis_NonlinearDifferentialEquations, Hunter_LaneEmden, KyciaFilipuk_LandEmden,  ISAAC, KyciaFilipuk_LandEmden}:
\begin{equation}
 \frac{d^{2}u(x)}{dx^{2}}+\frac{\alpha }{x} \frac{du(x)}{dx}+\delta u(x)^{p}=0,
 \label{Generlized_Lane-Emden}
\end{equation}
where $\alpha$ and $\delta$ are positive real constants and $p>1$ is natural number.

It is assumed that the perturbed equation of (\ref{Generlized_Lane-Emden}) is in the following form
\begin{equation}
\begin{array}{c}
p(x)\frac{d^{2}u(x)}{dx^{2}}+q(x)\frac{du(x)}{dx}+ r(x)u(x)+\delta u(x)^{p}=0, \\
p(x) = (1+a_{-1}x+a_{0}x^{2}+\ldots+a_{n}x^{n+2}), \\
q(x) = \left( \frac{\alpha}{x}+b_{-1}+b_{0}x+\ldots+b_{n}x^{n+1}\right),\\
r(x) = \left(c_{-1}\frac{1}{x}+c_{0}+\ldots+c_{n}x^{n}\right),
\end{array}
\label{General_equation}
\end{equation}
where $n > -1$ is natural number, $\alpha>0$ and $\delta \neq 0$ are real constants. The restriction imposed on the coefficients $\{a_{k}\}_{k=-1}^{n}$, $\{b_{k}\}_{k=-1}^{n}$  and $\{c_{k}\}_{k=-1}^{n}$ will be given below. The polynomials $p,q,r$ introduce into equations new (fixed) singularities in the complex plane and by suitable change of coordinates we can assume that the closes singularity to the origin is located at $x=1$, which results in
\begin{equation}
 p(x)=(1-x)\bar{p}(x)=(1-x)(1+\bar{a}_{-1}x+\ldots+\bar{a}_{n-1}x^{n}+\bar{a}_{n}x^{n+1}),
\end{equation}
where $\bar{p}(x)$ has zeros in the complex plane that are farther from the origin than $1$. The coefficients correspondence is given by $\bar{a}_{-1}=a_{-1}$, $a_{k}=\bar{a}_{k+1}-\bar{a}_{k}$ for $k\in \{0,\ldots,n-1\}$ and $a_{n}=-\bar{a}_{n}$.

\emph{It will be explained how to find the solution of boundary problem - find analytic solution(s) which connect these two fixed singularities at $x=0$ and $x=1$. The procedure that construct such kind of analytic solutions and therefore provide solution for this nonlinear (singular) boundary value problem will be the main topic of this paper.}

This problem can be used to construct solutions of (\ref{General_equation}) on the unit interval and match it for some $x>1$ with a solution to some other equation due to analyticity and therefore it is important in applications.

By deformation of the unit interval into some smooth curve in the complex plan we will get general problem of finding analytic solutions on such curve that connects two (not necessary closest) fixed singularities of transformed (\ref{General_equation}). The assumption on the smoothness of the curve is important as we do not want to introduce additional singularities along it.

The equation (\ref{General_equation}) is general enough to describe many problems that arise in mathematical physics. For example the equation
\begin{equation}
 (1-x^{2})\frac{d^{2}u(x)}{dx^{2}} + \left( \frac{\alpha}{x}+\beta x\right)\frac{du(x)}{dx} - \gamma u(x) + \delta u(x)^{p}=0,
 \label{Equation1}
\end{equation}
 appears in  \cite{metoda_dim=3, wyzsze_wymiary, wyzsze_wymiary_movable_singularities} for some special values of parameters as equation for self-similar (analytic) profiles of nonlinear wave equation in flat spacetime on the unit interval. In this paper the general theory for such kind of equations will be provided and (\ref{Equation1}) will be analysed as an example for real parameters: $\alpha >0$, $\beta$, $\gamma$ and $\delta \neq 0$.

The paper is organized as follows: In the next section the existence of some special singular solution at $x=0$, which is a crucial ingredient for global analytic solution on the interval $[0;1]$ will be given, and then  the local existence of the analytic solutions around singular point $x=0$ and $x=1$ will be proved. Then the next section presents the general discussion of the existence of global solutions on the unit interval. Finally, an example of application of the results of the paper will be presented on (\ref{Equation1}). In the Appendix some technical results are collected for the Reader's convenience.

\section{Local existence}
\label{Section_General_equation}

\subsection{Local solution at $x=0$}

For the existence of the global analytic solutions on the unit interval the following singular solution at $x=0$ will be of paramount importance 
\begin{equation}
 u_{\infty}(x)=b_{\infty}x^{a}, \quad b_{\infty}=\left(\frac{2[\alpha(p-1)-(p+1)]}{\delta(p-1)^{2}}\right)^{1/(p-1)}, \quad a=\frac{1}{1-p}.
 \label{u_inf_solution}
\end{equation}
Therefore first the class of equations (\ref{General_equation}) which has the solution (\ref{u_inf_solution}) is singled out by the following
\begin{Proposition}
\label{Proposition_existence_u_inf_solution}
 If the equation (\ref{General_equation}) allows the solution (\ref{u_inf_solution}) then the coefficients have to fulfil the following conditions
 \begin{equation}
  a(a-1) a_{k}+ab_{k}+c_{k}=0, \quad k \in \{-1,\ldots,n\},
  \label{matching_condition}
 \end{equation}
 where $a$ is defined in (\ref{u_inf_solution}).
\end{Proposition}
\begin{Proof}
 The proof is straightforward. Substituting (\ref{u_inf_solution}) into (\ref{General_equation}) and collecting terms of the same order in $x$ we obtain
 \begin{equation}
  [b_{\infty}(a(a-1)+\alpha a) + \delta b_{\infty}^{p}]x^{a-2}+\sum_{k=-1}^{n}b_{\infty}[a(a-1)a_{k}+a b_{k}+c_{k}]x^{a+k}=0.
 \end{equation}
 The term at $x^{a-2}$ vanishes, which is exactly the definition of $b_{\infty}$. In addition, if (\ref{u_inf_solution}) is the solution of (\ref{General_equation}) then all the other coefficients terms have to vanish, which gives the condition (\ref{matching_condition}). 
\end{Proof}
Hereafter we assume that the equation (\ref{General_equation}) fulfils (\ref{matching_condition}), and therefore (\ref{u_inf_solution}) is a solution.

The next proposition shows the existence of local analytic solutions at $x=0$ for (\ref{General_equation}) and their behaviour for large initial data $u(0)=c$.
\begin{Proposition}
\label{Proposition_General_Equation_spiral}
There exist analytic solutions of (\ref{General_equation}) at $x=0$ with initial data $u(0)=c$.

For large $c$ the asymptotic of this analytic solution of (\ref{General_equation}) at fixed $x_{0}$ within its radius of convergence is described by scaled asymptotics of the generalized Lane-Emden equation, namely (see also Eq. (\ref{Generalized_Lane-Emden_Asumptotics_for_Large_x}) in the Appendix),
 \begin{equation}
 u_{\pm}(x_{0},c)\approx b_{\infty} x_{0}^{-2/(p-1)}(1 \pm A_{0} c^{-\frac{\alpha+3+p(1-\alpha)}{4}} x_{0}^{\frac{\alpha+3+p(1-\alpha)}{2(p-1)}}\sin(\omega \ln(c^{\frac{p-1}{2}}x_{0})+\phi),
 \label{Generalized_Equationn_Asumptotics_for_Large_c}
\end{equation}
where $p \neq p_{Q}$, odd and $f(p,\alpha)<0$, where
 \begin{equation}
  p_{Q}:=\frac{\alpha+3}{\alpha-1},
  \label{critical_p}
 \end{equation}
\begin{equation}
 f(p,\alpha):=(-1+\alpha)^{2}+p^{2}(9-10\alpha + \alpha^{2})-2p(-3-6\alpha+\alpha^{2}),
 \label{Generalized_Lane_Emden_f(p,a)}
\end{equation}
and where
\begin{equation}
 \omega(p,\alpha)=i\frac{\sqrt{-f(p,\alpha)}}{2(p-1)}.
\end{equation}
Here $A_{0}$ and $\phi$ are constants.
\end{Proposition}
\begin{Proof}
 First, we have to prove that there exists local analytic solution at $x=0$. If we rewrite (\ref{General_equation}) in the form
 \begin{equation}
  \left\{
  \begin{array}{l}
   xu' = xv \\
   xv' = \frac{-1}{1+\sum_{k=-1}^{n}a_{k}x^{k+2}}\left(\alpha +v\sum_{k=-1}^{n}b_{k}x^{k+2}+ u\sum_{k=-1}^{n}c_{k+1}x^{k+1}+x\delta u^{p}\right)
  \end{array}
  \right.
 \end{equation}
 The second equation has the following form
 \begin{equation}
  xv' = -\alpha v + x g(x,u,v),
 \end{equation}
 where $g$ is analytic at $x=0$. Therefore, using Proposition \ref{Proposition1} from the Appendix, we conclude that there is analytic solution at $x=0$ expressible in a power series form, which first term is initial data $c=u(0)$. It can be found by formal procedure of substituting a formal power series $u(x)=\sum_{k=0}^{\infty}h_{k}x^{k}$ into (\ref{General_equation}) and obtaining recurrence for the coefficients $h_{k}$ of this series. Then above statement assures us that this formal series is convergent - it is a solution.
 
 For proving the asymptotic (\ref{Generalized_Equationn_Asumptotics_for_Large_c}), the equation (\ref{General_equation}) will be transformed into the generalized Lane-Emden equation (\ref{Generlized_Lane-Emden}) and the Proposition (\ref{Proposition_Lane_Emden_Asymptotics_For_Large_x}) of the Appendix will be used.
 We use the substitution from \cite{metoda_dim=3, wyzsze_wymiary}
 \begin{equation}
  u = c w, \qquad y=c^{\gamma}x,\qquad \gamma = (p-1)/2
  \label{General_equation_Lane-Emden_substitution}
 \end{equation}
  which transform the equation (\ref{General_equation}) into ($'=\frac{d}{dy}$)
 \begin{equation}
  w''+\frac{\alpha}{y}w'+\delta w^{p}=-\frac{1}{c^{p}}\left(w''\sum_{k=-1}^{n}a_{k}\frac{y}{c^{k\gamma+1}}+w'\sum_{k=-1}^{n}b_{k}\frac{y^{k+1}}{c^{k\gamma+1}}+\sum_{k=-1}^{n}c_{k}\frac{y^{k}}{c^{k\gamma+1}} \right).
  \label{General_equation_Lane-Emden_transformed}
 \end{equation}
  In the limit $c \rightarrow \infty$ the RHS of (\ref{General_equation_Lane-Emden_transformed}) vanishes and we are left with the Lane-Emden equation (\ref{Generlized_Lane-Emden}). To make the limiting process precise we start from analytic solution at a point $x_{0}$ within its circle of convergence of analytic solution. Then we perform $c \rightarrow \infty$ limit termwise which transform the analytic solution into the series solution (\ref{Generalized_Lane-Emden_Solution}) for the LE equation (\ref{Generlized_Lane-Emden}) - the solution of (\ref{General_equation_Lane-Emden_transformed}) with vanishing RHS and initial condition $w(0)=1$. Using large $y$ asymptotic of the LE equation (\ref{Generlized_Lane-Emden}) of the form (\ref{Generalized_Lane-Emden_Asumptotics_for_Large_x}) from the Appendix, and returning to the original variables $u$, $x$, we obtain exactly  (\ref{Generalized_Equationn_Asumptotics_for_Large_c}), as claimed.
\end{Proof}
The asymptotic (\ref{Generalized_Equationn_Asumptotics_for_Large_c}) is analogous to the Proposition \ref{Proposition_Lane_Emden_Asymptotics_For_Large_x} of teh Appendix, and results form the appearance of the LE equation (\ref{Generlized_Lane-Emden}) as a basic building block of (\ref{General_equation}). Loosely speaking, the Proposition \ref{Proposition_General_Equation_spiral} shows that when the phase plane $(u(x_{0});u'(x_{0}))$ fixed at $x_{0}$ is considered, then the analytic solution of (\ref{General_equation}) behaves as the spiral (\ref{Generalized_Equationn_Asumptotics_for_Large_c}) parametrized by $c$ that wraps around the limit point 
\begin{equation}
 P_{\infty} = \left(b_{\infty} x_{0}^{-2/(p-1)}, -\frac{2}{p-1}b_{\infty} x_{0}^{-2/(p-1)}\right). 
 \label{oscillation_point}
\end{equation}
The same behaviour shows special case of (\ref{General_equation}) and it was described in \cite{metoda_dim=3, wyzsze_wymiary}.

One can note that the Proof does not depend on the condition (\ref{matching_condition}).

An obvious Corollary for the Proposition \ref{Proposition_General_Equation_spiral} is 
\begin{Corollary}
 For large $u(0)=c$ value the analytic solutions have movable singularities located as for the analytic solutions of the LE equation (\ref{Generlized_Lane-Emden}), among other movable singularities.
\end{Corollary}
\begin{Proof}
 Since (\ref{General_equation_Lane-Emden_transformed}) for $c\rightarrow \infty$ goes to (\ref{Generlized_Lane-Emden}) therefore during this process the movable singularities of the analytic solution of (\ref{General_equation}) approach to the movable singularities of (\ref{Generlized_Lane-Emden}) described in \cite{KyciaFilipuk_LandEmden}, i.e., they are located symmetrically around the origin on the rays connecting the origin with all $p$ roots of $-1$.
\end{Proof}

This result ends our discussion of behaviour of local analytic solutions at $x=0$. We now pass on to the existence of the local solution at $x=1$.

\subsection{Local existence at $x=1$}
\label{SubSection_Local_Existence_x1}
In order to consider the local analytic solutions around $x=1$ the new variable $y=1-x$ is introduced. The equation (\ref{General_equation}) has now the form
\begin{equation}
\begin{array}{c}
 y\bar{p}(y)\frac{d^{2}u(y)}{dy^{2}}-q(y)\frac{du(y)}{dy}+ r(y)u(y)+\delta u(y)^{p}=0, \\
\bar{p}(y) = (A+A_{-1}y+A_{0}y^{2}+\ldots+A_{n}y^{n+2}), \\
q(y) = \left( \frac{\alpha}{1-y}+B_{-1}+B_{0}y+\ldots+B_{n}y^{n+1}\right),\\
r(y) = \left(c_{-1}\frac{1}{1-y}+C_{0}+\ldots+C_{n}x^{n}\right),
\end{array}
 \label{General_equation_y}
\end{equation}
where $A= 1+ \sum_{l=-1}^{n}a_{l}$ and $B_{-1}=b_{-1}-\sum_{l=0}^{n}b_{n}$.

On substituting formal powers series $u(y)=\sum_{l=0}^{\infty}d_{l}y^{l}$ and using the well-known Cauchy formula (\ref{power_of_power_series}) from the Appendix, the recurrence for $d_{l}$ coefficients is obtained
\begin{equation}
  \begin{array}{c}
   d_{0}=b, \\
   (l+1)(Al-B)d_{l+1}=f_{l+1}(d_{0},\ldots,d_{l}), \quad l>-1,
  \end{array}
  \label{recurrence_y0}
\end{equation}
where, as before, $A= 1+ \sum_{l=-1}^{n}a_{l}$, $B=\alpha + B_{-1}$ and $f_{l+1}$is a polynomial function and $b=u(y=0)$ is the initial data at $x=1$. Defining
\begin{equation}
 k=\frac{B}{A},
 \label{k_def}
\end{equation}
we have the following two cases
\begin{enumerate}
 \item {$k < 0$ or $k>0$, noninteger - there is infinite recurrence relations for the coefficients $d_{l}$;}
 \item {$k \in \mathbf{N}$ (resonance condition)- there is formal solution of the form $u(y)=d_{0}+ \ldots + d_{k}y^{k}+ d_{k+1}(d_{0},\ldots,d_{k},b)y^{k+1} +\ldots$, where the coefficients $\{d_{l}\}_{l=k+1}^{\infty}$ are derived from the recurrence relation and the first $k+1$ coefficients have fixed values that results from the set of polynomial equations
 \begin{equation}
  \left\{
  \begin{array}{c}
   f_{0}(d_{0})=d_{0}(C+d_{0}^{p-1})=0, \\
   f_{1}(d_{0},d_{1})=0, \\
   \ldots \\
   f_{k}(d_{0},\ldots, d_{k})=0,
  \end{array}
  \right.
 \end{equation}
 where $C =c_{-1}+C_{0}$. In principle such system can have complex-valued solutions.
 }
\end{enumerate}

The next step is to prove that these formal solutions are in fact solutions, i.e. they are convergent. It is immediate for finite-form solutions of the second case (when $b$ is selected in such a way that $d_{k+1}=0$). The structure of formal solutions is similar to those for (\ref{Equation1}) with some specially selected values of parameters, and the convergence of this special case was given in the Appendix of \cite{wyzsze_wymiary}. Therefore the proof for the convergence of formal power series of (\ref{General_equation_y}) is similar. The only difference is the multiplicity of coefficients. Therefore we have

\begin{Proposition}
\label{Proposition_General_Equation_existence_x1}
 Formal analytic solutions of (\ref{General_equation_y}) at $y=0$ obtained above are in fact the solutions, i.e., they are convergent in some neighbourhood of $y=0$.
\end{Proposition}

\begin{Proof}
 The idea of the proof is to transform (\ref{General_equation_y}) into the form required by the Proposition \ref{Proposition1} from \cite{lokalne_istnienie_praca} cited in the Appendix. The proof is similar to those of Appendix of \cite{wyzsze_wymiary} and therefore the sketch is only provided. 
 
 The proof depends on the value of $k$.
 
 For $k<0$, the equation (\ref{General_equation_y}) can be transformed into
 \begin{equation}
  \left\{
  \begin{array}{c}
   yu'=yv \\
   yv'= kv-Cu-\delta u^{p}+yf(y,u,v),
  \end{array}
  \right.
 \end{equation}
where $f(y,u,v)$ is some analytic function of arguments and $k$ is given by (\ref{k_def}). Introducing new variable $\tilde{v}$ via
\begin{equation}
 v=\tilde{v}+\frac{1}{kC}u+\frac{\delta}{k}u^{p},
\end{equation}
we get 
 \begin{equation}
  \left\{
  \begin{array}{c}
   yu'=yg(y,u,v) \\
   yv'= kv+y\tilde{f}(y,u,v),
  \end{array}
  \right.
 \end{equation}
 where again $g$ and $\tilde{f}$ are some analytic functions. Since $k<0$ we can apply aforementioned proposition which guarantees convergence of formal power series.
 
 For non-integer $\llcorner k \lrcorner < k < \ulcorner k \urcorner$ the general idea of the proof is as follows. First the new function $w(y)$ via
 \begin{equation}
  u(y)=d_{0}+\ldots + (d_{\llcorner k \lrcorner} + yw(y))y^{\llcorner k \lrcorner},
 \end{equation}
 is defined, where $\{d_{l}\}_{l=0}^{\llcorner k \lrcorner}$ coefficients are given by the recurrence  relation (\ref{recurrence_y0}). Then switching to the first order system of ODEs, removing constant terms and diagonalizing linear parts we arrive into the system required by the Proposition \ref{Proposition1}.

The last part of the proof is also valid for $k \in \mathbf{N}$.
\end{Proof}
The idea of the proof resembles the transformation to the normal form \cite{Wiggins, IlyaNormalForm}. As for $x=0$, one can note that the Proof does not depend on the condition (\ref{matching_condition}) and therefore it is a general result for (\ref{General_equation}).

In the next section general discussion on the existence of global analytic solutions on unit interval will be given.

\section{Global existence of analytic solutions}

In order to get global solution generalization of the method proposed in \cite{metoda_dim=3} and developed in \cite{wyzsze_wymiary} will be used. For such general class as (\ref{General_equation}) we present only steps that can be done and not a strict algorithm - these steps have to be adjusted to the concrete equation under consideration. The steps are as follows:

\begin{itemize}
 \item {Construct local analytic solutions at $x=0$ (existence proved in the Proposition \ref{Proposition_General_Equation_spiral} ) and prove that they can be extended towards $x=1$.}
 \item {Construct local analytic solutions at $x=1$ (existence proved in the Proposition \ref{Proposition_General_Equation_existence_x1} )  and prove that they can be extended towards $x=0$.}
 \item {Match theses local solutions at some point inside unit interval to obtain global solution.}
\end{itemize}

For proving that there is no singularities along the unit interval one can use some global methods, e.g., Lyapunov function (see \cite{Hartman, Wiggins} for general reference on Lyapunov functions). If this is proved then local analytic solution can be extended along real axis from $x=0$ towards $x=1$. Usually, the similar proof can be performed for local analytic solution that starts from $x=1$, and it can be shown that it can be extended to $x=0$ along real axis.

The proof of analytic continuation of local analytic solutions will be provided for a restricted class of (\ref{General_equation}), namely, we have

\begin{Proposition}
\label{Proposition_Lyapunov}
 For $\left( \frac{p'(x)}{2}-q(x)\right) < 0$, $p(x)>0$ in $x\in [0;1]$, a constant $r(x)=r$, and odd $p$, the local analytic solutions of (\ref{General_equation}) around $x=0$ and $x=1$ described in the previous section can be extended along the unit interval, i.e., there is no (movable) singularities of the extended solutions in $[0;1]$.
\end{Proposition}

The assumption $p(x)>0$ in the unit interval is not restrictive, since the whole equation can be multiplied by $-1$. There is $p(x) \neq 0$ in the unit interval since there is no (fixed) singularity of (\ref{General_equation}) in this interval by assumption.

\begin{Proof}
 The proof goes along the same line as in \cite{metoda_dim=3, wyzsze_wymiary} and base on the Lyapunov function 
 \begin{equation}
  H(x)=\frac{p(x)}{2}u'(x)^2+\frac{\delta}{p+1}u(x)^{p+1}+\frac{r}{2}u(x)^{2}.
 \end{equation}
 Its derivative on the solutions of (\ref{General_equation}) is
 \begin{equation}
  H'= \left( \frac{p'(x)}{2}-q(x)\right) u'^{2} <0.
 \end{equation}
 Since $H$ is monotone on solutions, therefore denoting  $V_{min}=min_{u} \left( \frac{\delta}{p+1}u^{p+1}+\frac{r}{2}u^{2} \right)$ ($p$ is odd and therefore global minimum exists) we have from $H(x)\leq H(0)$ the following estimate
 \begin{equation}
  \frac{p(x)}{2}u'(x)^2 < H(0) - V_{min},
 \end{equation}
 and therefore, since $p(x) >0$ in $[0;1]$ we get that $|u'(x)|$ is bounded and the analytic solution extended from $x=0$ towards $x=1$ is also bounded.
 
 For extension of the analytic solution at $x=1$ towards $0$ the following inequality will be used
 \begin{equation}
  \frac{-H'}{1+H}\leq - \frac{p'(x)-2q(x)}{p(x)}.
 \end{equation}
Integrating from $x=1$ to some $0<\bar{x} << 1$ and using the fact that the RHS is regular in this interval we get that $H$ is bounded and therefore $u$ and $u'$ are finite as well. 
\end{Proof}

\begin{figure}[htb!]
\centering
 \includegraphics[width = 0.5\textwidth]{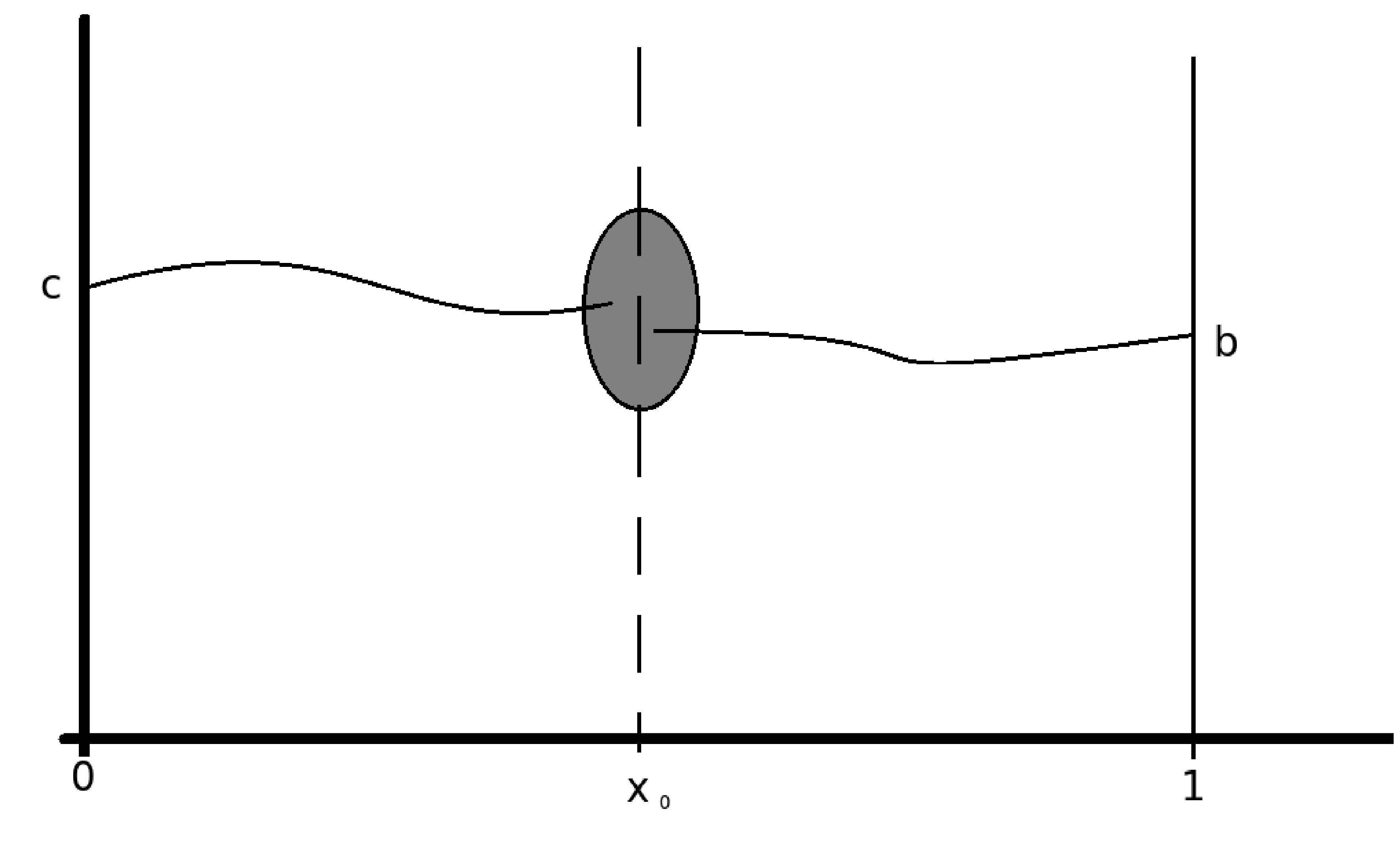}
 \caption{Matching of two local analytic solutions extended to $x_{0}$ point.}
 \label{Figure_matching1}
\end{figure}

\begin{figure}[htb!]
\centering
 \includegraphics[width = 0.5\textwidth]{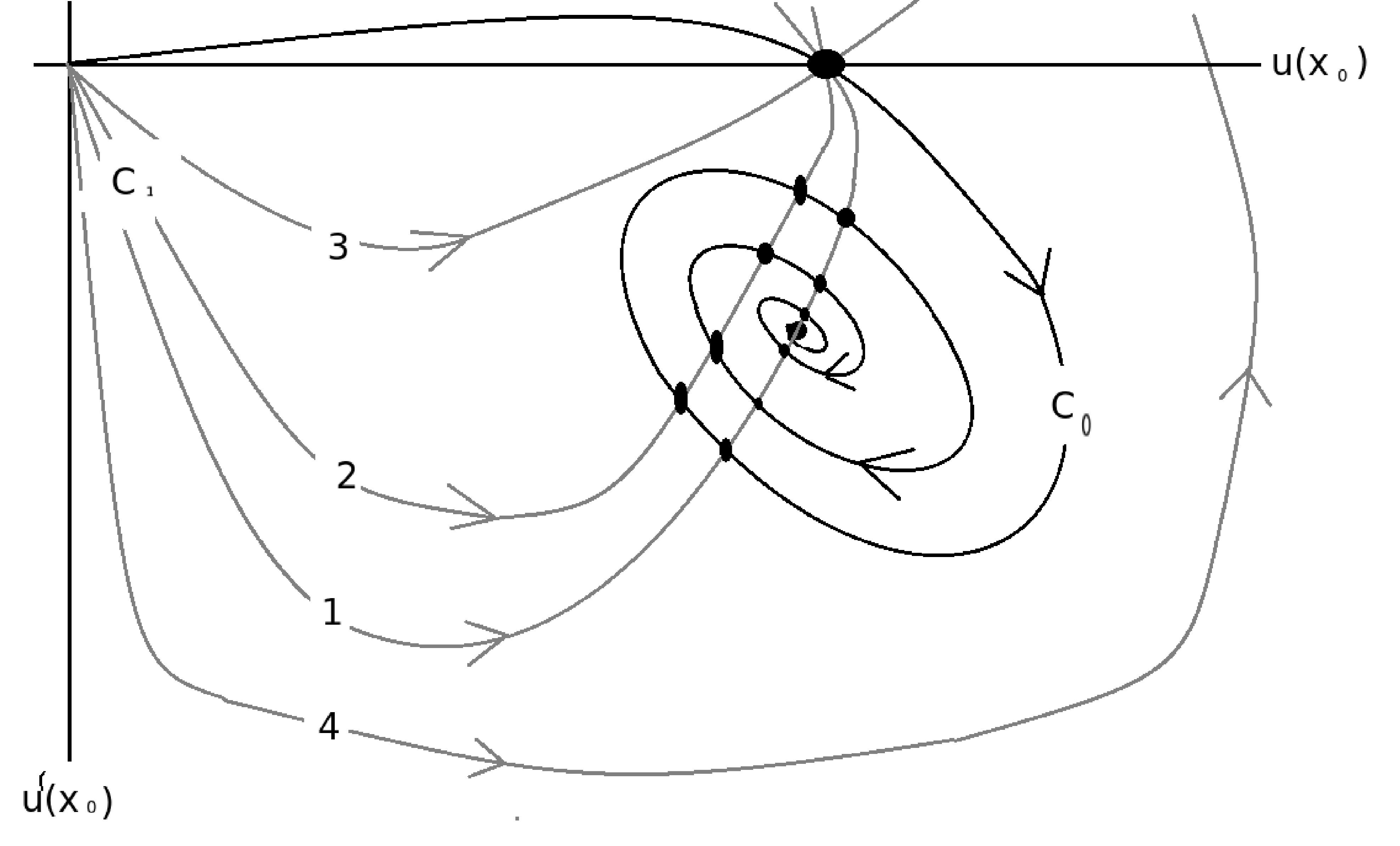}
 \caption{Matching of two local analytic solutions extended to $x_{0}$ point. The cases $1$,$2$,$3$ belong to the same type of the equation with different values of coefficients. The case $4$, topologically different from the previous ones corresponds to the different equation. Wrapping point of $C_{0}$ curve moves when the values of the equation's coefficients are changing.}
 \label{Figure_matching2}
\end{figure}

In the last step analytic solutions extended from both endpoints has to be matched at some intermediate point $x_{0}\in [0;1]$, as it is presented in Fig. \ref{Figure_matching1}.
The matching can be intuitively visualized when we analyse situation on the phase plane at $x_{0}$. Fig. \ref{Figure_matching2} presents the situation.
If we will vary $c$ initial data of the solution starting form $x=0$ then, as Proposition \ref{Proposition_General_Equation_spiral} states, we obtain the  spiral that wraps around $P_{\infty}$, which is $C_{0}$ curve.  The second curve $C_{1}$ is the image of the solution extended form $x=1$, and it is parametrized by $u(1)=b$ - the initial data at this endpoint. Each intersection of these curves is exactly $C^{1}$ matching condition, moreover, it is analytical matching, because we match two analytic solutions.

Going along $C_{1}$ curve we start from $b=0$ - the trivial solution. Increasing $b$ we can get a few possibilities for the number of intersections: 
\begin{enumerate}
 \item {In the first case the $C_{1}$ curve passes exactly through the wrapping point of the spiral - $P_{\infty}$. In this case we obtain countable family of analytic solutions on unit interval (case 1 in the figure). For this case the existence $u_{\infty}$ solution, i.e., the condition (\ref{matching_condition}), is required. It is obvious to explain why it is the case if we use continuous dependence of solutions on initial conditions with initial data that are slightly shifted form singular points $x=0$ and $x=1$, which however match these analytic solutions for some initial data specified at these singular points. The $x=0$ solutions 'wraps around' (\ref{u_inf_solution}). The second solution from $x=1$ stays close to this point if we vary initial data $b$ around $b_{\infty}$ value, therefore, also around (\ref{u_inf_solution}).}
 \item {The second curve can cross $C_{0}$ in finite number of points - finite number of global solutions (case 2 in the figure). In this case the condition (\ref{matching_condition}) is slightly perturbed.}
 \item { If there are some special solutions, like constant ones for the case of resonances for  $x=1$ solution (see also next section for example), then the second curve can cross $C_{0}$ in one point which corresponds to this special solution (case 3 in the figure).}
\item {The second curve can also have only one point $(0;0)$ (which corresponds to the trivial solution) with the $C_{0}$ curve (case 4 in the figure).}
\end{enumerate}

The statements of this qualitative discussion will now be explained on the example in the following section.

\section{Example}
\label{Section_Example}
In this section a simple example (\ref{Equation1}) that illustrate above results and discussion is presented. It is a generalization of the equation analysed in \cite{metoda_dim=3, wyzsze_wymiary, wyzsze_wymiary_movable_singularities}. The other example of the class (\ref{General_equation}) is discussed in \cite{Example}.

It is assumed that $\alpha >0$, $\beta$, $\gamma$ and $\delta \neq 0$ are real parameters for (\ref{Equation1}).

The condition for the existence of the solution (\ref{u_inf_solution}) is of the form 
\begin{equation}
 \gamma = \frac{2}{1-p}\left( \frac{p+1}{p-1}+\beta\right),
 \label{Equation1_u_inf_existence_condition}
\end{equation}
and it is assumed that it is hereafter fulfilled. In addition, there exists the constant solution
\begin{equation}
 u_{0}=\left( \frac{\gamma}{\delta}\right)^{1/(p-1)}:=b_{0},
 \label{u_0_solution}
\end{equation}
providing that $\gamma >0$ ($\delta >0$ was assumed previously).

We start from analysing local power series solution at $x=0$, which the Proposition \ref{Proposition_General_Equation_spiral} assures to exist. On substituting formal power series $u(x)=\sum_{k=0}^{\infty}a_{k}x^{k}$ into (\ref{Equation1}) we get unique recurrence for the coefficients
\begin{equation}
 a_{0}=c=\text{arbitrary}, \quad a_{1}=0,\quad a_{k+2}=\frac{[(k(k-1-\beta)+\gamma]a_{k}-\delta c_{k}}{(k+2)(k+1+\alpha)},
 \label{Equation1_RecurrenceRealtion_at_x=0}
\end{equation}
where $k \ge 0$ and $c_{k}$ can be computed from $b_{k}$ ones using (\ref{power_of_power_series}) formula. This power series can be analytically continued from $x=0$ to the vicinity of $x=1$ along real line without encounter any singularities. It can be proved by repeating the Proof of  the Proposition \ref{Proposition_Lyapunov} with additional assumptions, that $p$ is odd, and $\delta >0$. To this end we use the following Lyapunov function
\begin{equation}
 H(u,u';x)=(1-x^{2})\frac{u'^{2}}{2}-\gamma \frac{u^{2}}{2}+\frac{\delta}{p+1}u^{p+1}.
 \label{Equation1_Lyapunov_function}
\end{equation}
We can shift this function by a constant (the minimum of polynomial in $u$ consisting of the last two terms) to make it positively defined.
Taking derivative and replacing $u''$ with the help of (\ref{Equation1}) we obtain
\begin{equation}
 H'= - \left( (1+\beta)x +\frac{\alpha}{x}\right) u'^{2}.
\end{equation}
the term in the bracket is positive on the unit interval when $\beta>-1$ or $\beta=0$ or when $\beta<-1$ and $\alpha+\beta > -1$. Assuming this we have that $H'<0$, therefore, $H(0) > H(x)$ for $0<x<1$. From this we obtain immediately, as in \cite{metoda_dim=3}, \cite{wyzsze_wymiary} that solution which is regular at $x=0$ stays finite along the unit interval.

The solution at $x=1$ can be obtained by introducing the variable $y=1-x$ into (\ref{Equation1}), to obtain ($'=d/dy$)
\begin{equation}
 y(2-y)u''-\left( \frac{\alpha}{1-y}+\beta (1-y) \right) u' -\gamma u+\delta u^{p}=0. 
\end{equation}
Introducing a power series ansatz $u(y)=\sum_{k=0}^{\infty}b_{k}y^{k}$ we obtain the following recurrence for the coefficients
\begin{equation}
\begin{array}{l}
 b_{0}=b=\text{arbitrary}, \quad b_{1}=\frac{\gamma b_{0}-\delta c_{0}}{-(\alpha + \beta)},\\
 b_{k+1}=\frac{k(3(k-1)-2\beta)b_{k}+\gamma b_{k}-\delta c_{k}+(k-1)(\beta+2-k)b_{k-1}-\gamma b_{k-1}+\delta c_{k-1}}{(k+1)(2k-\alpha-\beta)}),
\end{array}
\label{Equation1_RecurrenceRealtion_at_x=1}
\end{equation}
where $c_{k}$ can be computed from $b_{k}$ ones employing (\ref{power_of_power_series}) formula. When $\alpha + \beta =2k$, where  $k \ge 0$ is a natural number we obtain the resonance condition, which leads, as in \cite{wyzsze_wymiary}, to the special class of the solutions described in Section \ref{SubSection_Local_Existence_x1}. We will not discuss this class as it is analogous to the analysis in \cite{wyzsze_wymiary}. Therefore, we assume that $-1< \alpha + \beta < 0$.

The convergence proof for the solution with the coefficients (\ref{Equation1_RecurrenceRealtion_at_x=1}) is given by the Proposition \ref{Proposition_General_Equation_existence_x1}, however it is instructive to provide special case of this proof. To prove that this power series has nonzero radius of convergence we rewrite it as a system of first order ODEs
\begin{equation}
\left\{
\begin{array}{l}
 yu'=yv \\
 yv' = \frac{1}{(2-y)(1-y)} ((\alpha +\beta (1+y^{2}-2y))v +\gamma(1-y)u-\delta (1-y)u^{p}.
\end{array}
 \right.
 \label{Equation1_ODE_system_at_x=1}
\end{equation}
We are allowed to use the Proposition \ref{Proposition1} if we remove the terms of $u$ and $u^{p}$ with constant coefficients. This can be done using the following change of variables
\begin{equation}
 v = \bar{v} - \frac{\gamma}{\alpha+\beta}u +\frac{\delta}{\alpha+\beta}u^{p},
\end{equation}
which transforms (\ref{Equation1_ODE_system_at_x=1}) into 
\begin{equation}
\left\{
\begin{array}{l}
 yu'=yf(y,u,v) \\
 yv' = 2(\alpha+\beta) + yg(y,u,v),
\end{array}
 \right.
\end{equation}
where $f$ and $g$ are functions analytic in all of its variables around $y=0$. From the fact that $\alpha+\beta <0$ and Proposition \ref{Proposition1} we conclude that the solution (\ref{Equation1_RecurrenceRealtion_at_x=1}) is convergent in some neighbourhood of $y=0$ ($x=1$).

When $2k<\alpha+\beta < 2(k+1)$ or when resonance condition occurs, namely $\alpha+\beta=2k$, for natural $k$, then the general proof is a simple modification of those from Appendix of \cite{wyzsze_wymiary}.

In next step we show, that the solution form $x=1$ can be extend along the real line to $x=0$. We use, as in \cite{metoda_dim=3} and \cite{wyzsze_wymiary}, the following easy to check bound
\begin{equation}
 \frac{-H'}{H+1} < \frac{2\alpha}{x},
\end{equation}
from which it occurs, by integrating both sides, that $u$ and $u'$ are bounded inside the unit interval if we start from solution (\ref{Equation1_RecurrenceRealtion_at_x=1}) bounded at $x=1$.

These results assures us that these two local solutions can be matched at some $x_{0}$ inside the unit interval. This matching fits exactly the general scheme presented in the previous section: On the phase plane at $x_{0}$ the solution form $x=0$ gives a spiral which wraps around (\ref{u_inf_solution}) solution. The second curve for solution at $x=1$ crosses this spiral, such that, every intersection gives a global solution.  We can note that there will be at least two such intersections, which corresponds to the $u=0$ and $u_{0}$ solution.

\begin{figure}[htb!]
\centering
 \includegraphics[width = 0.5\textwidth]{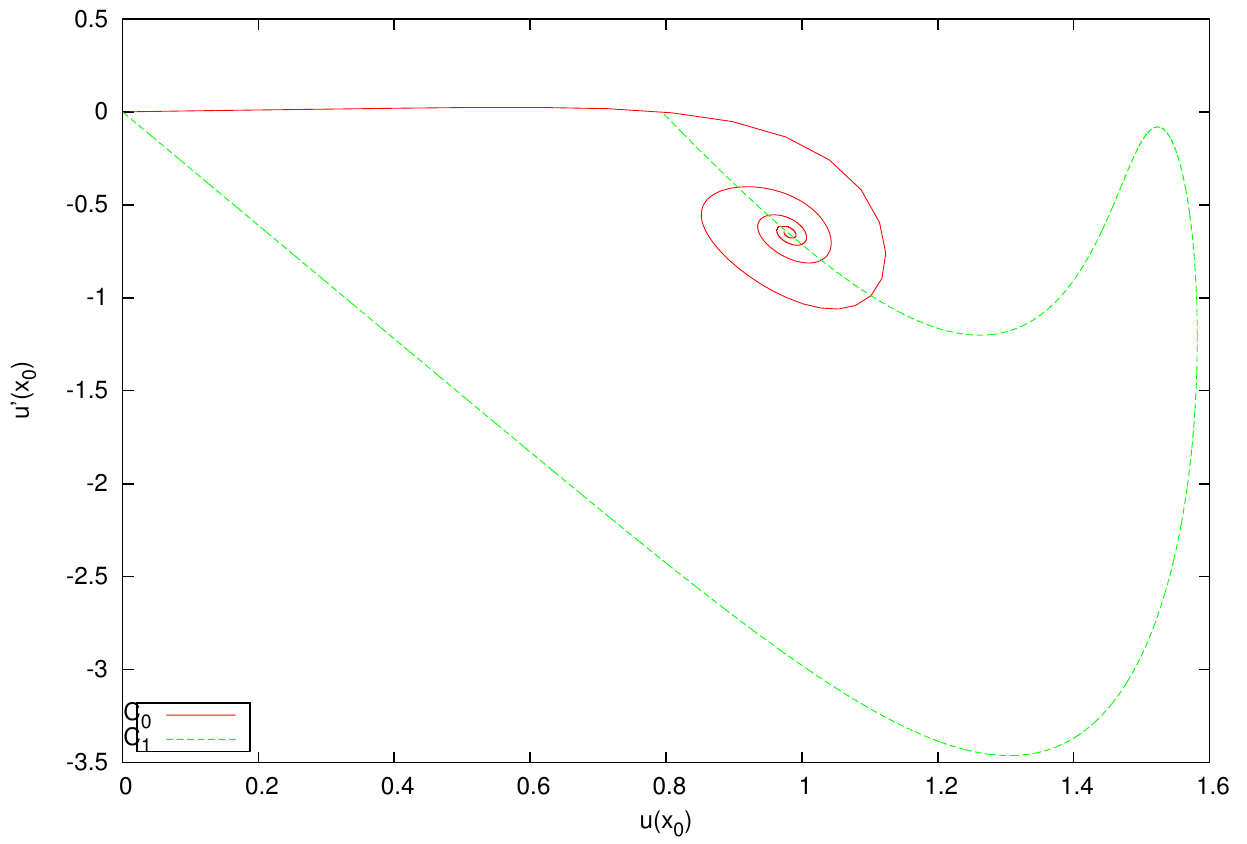}
 \caption{Plot for $\alpha =2$, $\beta=-25/12$, $\gamma=1/4$, $\delta =1$. In this case matching point $x_{0}=0.5$, and wrapping point $P_{\infty}=(\frac{\sqrt{2}}{\sqrt[3]{3}}; -\frac{2 \sqrt{2}}{3 \sqrt[3]{3}})$. The curve starts form $(0;0)$ point passes through $P_{\infty}$ and integration was terminated at $(b_{0};0)$ point.}
 \label{Figure_Equation1_countable_family}
\end{figure}

\begin{figure}[htb!]
\centering
 \includegraphics[width = 0.5\textwidth]{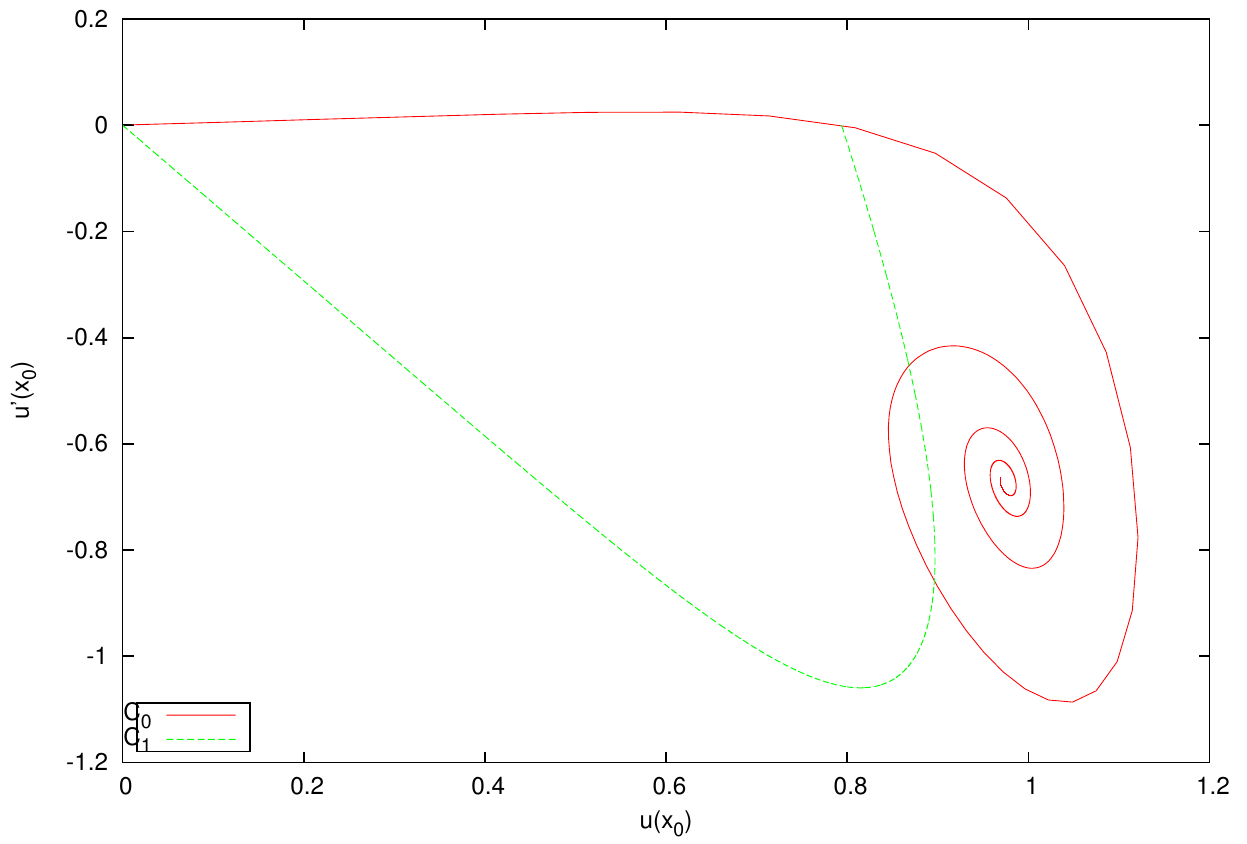}
 \caption{Plot for $\alpha =2$, $\beta=-29/12$, $\gamma=1/4$, $\delta =1$. In this case matching point $x_{0}=0.5$, and wrapping point $P_{\infty}=(\frac{\sqrt{2}}{\sqrt[3]{3}}; -\frac{2 \sqrt{2}}{3 \sqrt[3]{3}})$. The curve starts form $(0;0)$ point but not passes through $P_{\infty}$ as condition (\ref{Equation1_u_inf_existence_condition}) is not fulfilled. Integration was terminated again at $(b_{0};0)$ point.}
 \label{Figure_Equation1_finite_family}
\end{figure}

\begin{figure}[htb!]
\centering
 \includegraphics[width = 0.5\textwidth]{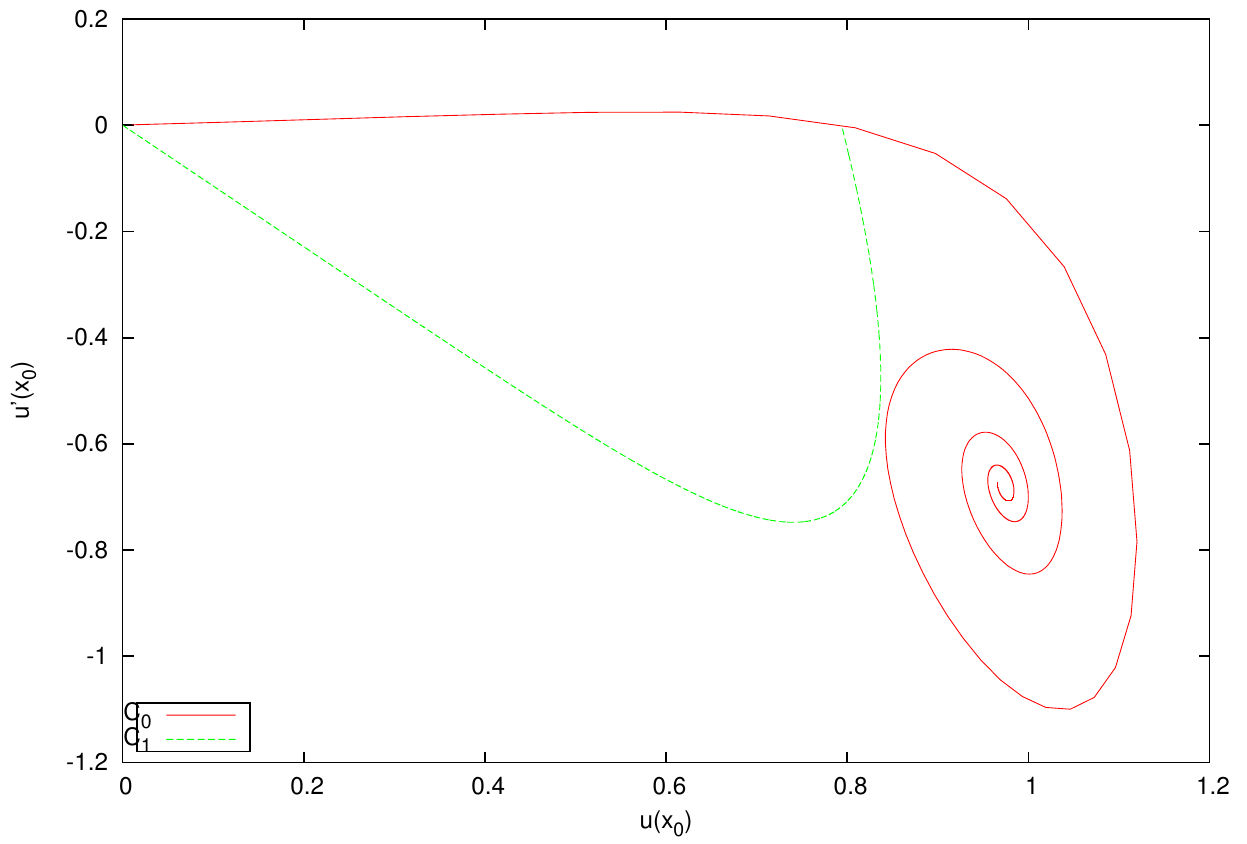}
 \caption{Plot for $\alpha =2$, $\beta=-31/12$, $\gamma=1/4$, $\delta =1$. In this case matching point $x_{0}=0.5$, and wrapping point $P_{\infty}=(\frac{\sqrt{2}}{\sqrt[3]{3}}; -\frac{2 \sqrt{2}}{3 \sqrt[3]{3}})$. The curve starts form $(0;0)$ point but not passes through $P_{\infty}$ as condition (\ref{Equation1_u_inf_existence_condition}) is not fulfilled. Integration was terminated again at $(b_{0};0)$ point. Now we only get two global solutions - the trivial one and $u_{0}$.}
 \label{Figure_Equation1_two_solutions}
\end{figure}

The Figure \ref{Figure_Equation1_countable_family} presents the case with countable many global solutions. For this case the condition (\ref{Equation1_u_inf_existence_condition}) is fulfilled. If we now slightly perturb the equations coefficients that this condition is not fulfilled then we get finite number of intersections/global solutions, and in the end we are left only with two constant solution. This is presented in Figs. \ref{Figure_Equation1_finite_family} and \ref{Figure_Equation1_two_solutions}.

This example shows that for the specific case presented in \cite{metoda_dim=3, wyzsze_wymiary} the coefficients appears perfectly to fulfil the matching condition (\ref{Equation1_u_inf_existence_condition}).

\section{Conclusions}
General method of matching using singular solution at one of the fixed singularity of the second order differential equation that generalize the Lane-Emden equation was presented. The method is constructive and shows under what circumstances countable family of analytic solutions that interpolate between fixed singularities of the equation exists. This is significant generalization of the idea presented from \cite{metoda_dim=3} and \cite{wyzsze_wymiary}.

This theory was used to describe matching problem for the differential equation generalizing the equations for self-similar profiles of nonlinear wave equations with power-type nonlinearity. This example shed a new light on the special case presented in \cite{wyzsze_wymiary}.

The boundary value problem and equations of the type presented in the paper are very common in applications since they arise form radial part of the Laplace operator and power-type nonlinearity, therefore, these results are useful in applied sciences and technical problems.

\section*{Acknowledgments}

This research was supported by the GACR grant 17-19437S, and the grant MUNI/A/1138/2017 of Masaryk University. 
We also thank the PHAROS COST Action (CA16214) for partial support.

\appendix

\section{Preliminary facts}
This Appendix collects some results useful in the course of the paper and scattered in the literature. They were collected for the Reader's convenience.

In the paper to get the existence of the local analytic solutions around finite fixed singularities Proposition 1 from \cite{lokalne_istnienie_praca} is employed, which reads

\begin{Proposition} \cite{lokalne_istnienie_praca}
\label{Proposition1}
Consider a system of differential equations for $i+j$ functions $u=(u_{1},\ldots,u_{i})$ and $v=(v_{1},\ldots,v_{j})$,
\begin{equation}
t\frac{du_{l}}{dt}=t^{\mu_{l}}f_{l}(t,u,v), \qquad t\frac{dv_{l}}{dt}=-\lambda_{l}v_{l}+t^{\nu_{l}}g_{l}(t,u,v),
\label{local_existence_system_preposition_1}
\end{equation}
with constants $\lambda_{l} >0$ and integers $\mu_{l},\nu_{l}\geq 1$ and let $U$ be an open subset of $R^{n}$ such that the functions $f$ and $g$ are analytic in a neighborhood of $t=0$, $u=c$, $v=0$ for all $c \in U$. Then there exists an $i$-parameter family of solutions of the system (\ref{local_existence_system_preposition_1}) such that 
\begin{equation}
 u_{l}(t)=c_{l}+O(t^{\mu_{l}}), \qquad v_{l}(t)=O(t^{\nu_{l}}),
\label{local_existence_solution_preposition_1}
\end{equation}
where $u_{l}(t)$ and $v_{l}(t)$ are defined for $c \in U$, $|t|<t_{0}(c)$ and are analytic in $t$ and $c$.
\end{Proposition}

In derivation of the recurrence for the equations with power types nonlinearities it is convenient to use the well-known Cauchy product \cite{Gradhsteyn_Ryzhyk}
\begin{equation}
\begin{array}{l}
\left( \sum_{l=0}^{\infty} a_{l}(x-x_{0})^{l} \right )^{p}=\sum_{l=0}^{\infty} c_{l}(x-x_{0})^{l}, \\ 
\\
c_{0}=a_{0}^{p}, \qquad c_{m}=\frac{1}{ma_{0}}\sum_{l=1}^{m}(lp-m+l)a_{l}c_{m-l},
\end{array}
\label{power_of_power_series} 
\end{equation}
for $m>0$, and where $a_{0}=c$ is a free parameter. The Cauchy product allows us to deal with nonlinear term in simple manner.

Some facts from the theory of (\ref{Generlized_Lane-Emden}) are presented following \cite{ISAAC}, \cite{KyciaFilipuk_LandEmden}. The equation (\ref{Generlized_Lane-Emden}) has an analytic solution around the fixed singularity $x=0$
\begin{Proposition} \cite{ISAAC}
 The equation (\ref{Generlized_Lane-Emden}) has local analytic solution around fixed singularity at $x=0$ of the form
 \begin{equation}
  a_{0}=\text{initial data}, \quad a_{1}=0, \quad a_{k+2}=\frac{\delta c_{k}}{(k+2)(k+1+\alpha)} \quad k \geq 0;
  \label{Generalized_Lane-Emden_Solution}
 \end{equation}
 The $\{c_{k}\}_{k=0}^{\infty}$ are coefficients derived form $\{a_{k}\}_{k=0}^{\infty}$ using (\ref{power_of_power_series}).
\end{Proposition}
The proof can be found in \cite{ISAAC}.

The equation (\ref{Generlized_Lane-Emden}) has also the singular solution (\ref{u_inf_solution}). It is singular at the origin and it somehow, attracts solutions that vanish at infinity. The precise meaning is given by the following Proposition which is generalization of the results form \cite{metoda_dim=3}, \cite{wyzsze_wymiary}.
\begin{Proposition} \cite{Hunter_LaneEmden, KyciaFilipuk_LandEmden}
\label{Proposition_Lane_Emden_Asymptotics_For_Large_x}
 For $p \neq p_{Q}$, odd and $f(p,\alpha)<0$, where
 \begin{equation}
  p_{Q}:=\frac{\alpha+3}{\alpha-1},
 \end{equation}
\begin{equation}
 f(p,\alpha):=(-1+\alpha)^{2}+p^{2}(9-10\alpha + \alpha^{2})-2p(-3-6\alpha+\alpha^{2})
\end{equation}
the asymptotic of analytic solution (\ref{Generalized_Lane-Emden_Solution}) with $u(0)=a_{0}=1$ is given in the following form
\begin{equation}
 u_{\pm}(x)\approx b_{\infty} x^{-2/(p-1)}(1 \pm A_{0} x^{\frac{\alpha+3+p(1-\alpha)}{2(p-1)}}\sin(\omega \ln(x)+\phi)
 \label{Generalized_Lane-Emden_Asumptotics_for_Large_x}
\end{equation}
for large $x$, where
\begin{equation}
 \omega(p,\alpha)=i\frac{\sqrt{-f(p,\alpha)}}{2(p-1)}
\end{equation}
and $A_{0}$ and $\phi$ are constants.
\end{Proposition}
The proof for $\alpha = 2$ and $\delta = 1$ can be found in \cite{Hunter_LaneEmden} and general proof in \cite{KyciaFilipuk_LandEmden}.




\end{document}